\documentclass[11pt]{amsart}

\usepackage{amsmath,amsthm, amssymb}
\usepackage[T1]{fontenc}
\usepackage{latexsym}
\usepackage{subfigure, epsfig, epsf}
\usepackage{color,fancyhdr,lastpage,graphicx}
\usepackage{xspace}
\usepackage{wasysym}
\usepackage{setspace}
\usepackage{bbm}
\usepackage{stmaryrd}

\usepackage{pstricks}
\usepackage{pst-node}
\usepackage{pst-tree}

    \newcommand{\QQ}{\mathbb{Q}}  
\newcommand{\RR}{\mathbb{R}}

 \newcommand{\mcC}{\mathcal{C}}  \newcommand{\mcF}{\mathcal{F}}

  \newcommand{\mcL}{\mathcal{L}}

\newcommand{\EE}{\mathbb{E}} \newcommand{\PP}{\mathbb{P}}

\newcommand{\ep}{\epsilon}
\renewcommand{\leq}{\leqslant}
\renewcommand{\geq}{\geqslant}

\newcommand{\al}{\alpha}

\newcommand{\Vol}{\textsc{Vol}}

\newcommand{\wrt}{with respect to }

\newcommand{\ssk}{\smallskip}

\newtheorem{thm}{\hspace{-0.15cm}  {\sc Theorem} }

\newtheorem{rems}[thm]{{\sc Remarks}}

\numberwithin{equation}{section} 

\newenvironment{Dem}{%
    \begin{list}{\hspace{0.5cm}{\sc Proof --}}{%
        \setlength{\topsep}{0pt}%
        \setlength{\leftmargin}{0pt}%
        \setlength{\rightmargin}{0pt}%
        \setlength{\listparindent}{0pt}%
        \setlength{\itemindent}{0pt}%
        \setlength{\parsep}{0pt}%
        \addtolength{\leftmargin}{20pt}%
        \addtolength{\rightmargin}{0pt}%
    } \item }{\hfill{\space $\rhd$}\end{list}\smallskip}

\title{Large deviation principle for bridges of degenerate diffusion processes}
\date{\today}
\author{I. Bailleul}
\address{IRMAR, 263 Avenue du General Leclerc, 35042 RENNES, France}
\email{ismael.bailleul@univ-rennes1.fr}
\thanks{This research was partially supported by an ANR grant "Retour post-doctorant".}

\begin{document}

\maketitle

\begin{abstract}
We prove that bridges of subelliptic diffusions on a compact manifold, with distinct ends, satisfy a large deviation principle in the space of H\"older continuous functions, with a good rate function, when the travel time tends to $0$. This leads to the identification of the deterministic first order asymptotics of the distribution of the bridge under generic conditions on the endpoints of the bridge.
\end{abstract}

\section{Introduction}

Let $M$ be a compact, connect and oriented $m$-dimensonal smooth manifold, and $V_1,\dots,V_\ell$ be smooth vector fields on $M$, whose Lie algebra has maximal rank everywhere. Given another vector field $V$ on $M$, set
\begin{equation}
\label{EqDefnL}
\mcL = \frac{1}{2}\sum_{i=1}^\ell V_i^2 +V.
\end{equation}
The semi-group associated with $\mcL$ has a smooth positive fundamental solution $p_t(z,z')$ with respect to any smooth volume measure $\Vol$ on $M$. Given $x,y$ in $M$ denote by $\Omega^{x,y}$ the set of continuous paths $\omega : [0,1]\rightarrow M$ with $\omega_0=x$ and $\omega_1=y$. For $\ep>0$, we define uniquely a probability measure $\PP^{x,y}_\epsilon$ on $\Omega^{x,y}$ defining $\PP^{x,y}_\epsilon\big(\omega_{t_1}\in A_1,\dots, \omega_{t_k}\in A_k\big) $ for all $k\geq 1,\, 0<t_1<\cdots<t_k<1$ and any Borel sets $A_1,\dots,A_k$ of $M$, by the formula
\begin{equation}
\label{EqConditonnedProbability}
\frac{1}{p_\epsilon(x,y)}\int \left(\prod_{j=1}^k \big(p_{t_j\ep - t_{j-1}\ep}(x_{j-1},x_j){\bf 1}_{A_j}(x_j)\right)\;p_{\epsilon-t_k\ep}(x_k,y)\,\Vol(dx_1)\cdots \Vol(dx_k)
\end{equation}
where $t_0=0$ and $x_0=x$. This formula describes the law of the diffusion process associated with $\epsilon\mcL$, conditionned on having position $x$ at time $0$ and position $y$ at time $1$. By Whitney's embedding theorem, there is no loss of generality in supposing that $M$ is a submanifold of an ambiant Euclidean space $\big(\RR^d,\|\cdot\|_d\big)$.

\ssk

Write $H^1_0$ for the set of $\RR^\ell$-valued paths $h$ over the time interval $[0,1]$, with starting point $0$; its $H^1$-norm is denoted by $\|h\|$. Given $h\in H^1_0$, we define a path $\gamma^h$ by solving the differential equation
\begin{equation}
\label{EqControlledODE}
\dot \gamma^h_t = \sum_{i=1}^\ell V_i\big(\gamma^h_t\big) \dot h^i_t,
\end{equation} 
for $0\leq t\leq 1$, given any specified starting point. The Lie bracket condition ensures that one defines a metric topology identical to the manifold topology setting for any pair of points $(a,b)$ in $M$
$$
d(a,b) = \inf \int_0^1 |\dot h_s|_\ell ds
$$
where the infimum is over the non-empty set of $H^1_0$-controls $h$ such that $\gamma^h_0=a$ and $\gamma^h_1=b$. It is called the sub-Riemannian distance associated with $\mcL$. The notation $|\cdot|_\ell$ stands here for the Euclidean norm on $\RR^\ell$. We define an $[0,\infty]$-valued function $J$ on $\Omega^{x,y}$ setting
\begin{equation}
\label{EqDefnJ}
J(\gamma) = \frac{1}{2}\Big(\inf\big\{\|h\|^2\,;\,\gamma^h=\gamma\big\} - d(x,y)^2\Big),
\end{equation}
with the convention $\inf\emptyset = \infty$. The above infimum is called the {\bf energy of the path $\gamma$}, classically denoted by $2I(\gamma)$.  Given any $0<\alpha<1$, denote by $\|x\|_\alpha$ the $\alpha$-H\"older norm of a path $x$ from $[0,1]$ to the ambiant space $\RR^d$. Write $\mcC^\alpha_{x,y}\big([0,1],M\big)$ for the set of all $M$-valued paths with finite $\alpha$-H\"older norm, with endpoints $x$ and $y$;  it is equipped with the topology associated with $\|\cdot\|_\alpha$.

\begin{thm}[Large deviation principle for bridges of degenerate diffusion processes]
\label{MainLDPThm}
~
\begin{enumerate}
   \item[(i)] Given any $\frac{1}{3}<\alpha<\frac{1}{2}$, the probabilities $\PP^{x,y}_\epsilon$ are supported on $\mcC^\alpha_{x,y}\big([0,1],M\big)$.   
   \item[(ii)] The family $\big(\PP^{x,y}_\epsilon\big)_{0<\ep<1}$ satisfies a large deviation principle in $\mcC^\alpha_{x,y}\big([0,1],M\big)$, with good rate function $J$.
\end{enumerate}
\end{thm}

\begin{rems}
\begin{enumerate}
   \item The above definition of the space $\mcC^\alpha_{x,y}\big([0,1],M\big)$ rests on the ambiant Euclidean metric. It is straigtforward to see that it coincides with the set of $M$-valued paths which are $\alpha$-H\"older for any choice of Riemannian metric on $M$, so $\mcC^\alpha_{x,y}\big([0,1],M\big)$ is intrinsically defined.  \vspace{0.2cm}
   \item Inahama proved in \cite{Inahama} a similar result under a stronger ellipticity condition. His analysis rests on the dynamic description of the diffusion associated with $\mcL$, given by the stochastic differential equation $dx_t = V(x_t)dt + \sum_{i=1}^\ell V_i(x_t){\circ dB^i_t}$, or rather on its rough path counterpart. By using quasi-sure analysis, he is able to lift the measures $\PP^{x,y}_\ep$ to some measures ${\bf P}^{x,y}_\ep$ on the space of geometric rough paths, which requires the quasi-sure existence of the Brownian rough path. The large deviation principle for $\PP^{x,y}_\ep$ is then obtained as a consequence of a subtle large deviation principle for ${\bf P}^{x,y}_\ep$, as the Ito-Lyons map is continuous. Our proof is more analytic, in that its essential ingredients are the heat kernel estimates of L\'eandre and Sanchez-Calle. We also use the machinery of rough paths as a convenient tool for proving the exponential tightness of the family of probability measures $\big(\PP^{x,y}_\epsilon\big)_{0<\epsilon\leq 1}$ on $\mcC^\alpha_{x,y}\big([0,1],M\big)$. 
   
   As a matter of fact, the proof below seems to be the first explicit proof of the above large deviation principle under the general Lie bracket condition for $\mcL$. It seems possible however to trace back the large deviation upper bound to some works of Gao \cite{Gao} and Gao and Ren \cite{GaoRen} on large deviation principles for stochastic flows in the framework of capacities on Wiener space. They prove in these works a Freidlin-Wentzell estimate/large deviation principle for $(r,p)$-capacities on Wiener space. Denote by $X^\ep$ the solution to the stochastic differential equation $dX^\ep_t = \ep\,V(X^\ep_t)dt + \ep^{1/2}\,V_i(X^\ep_t)\,{\circ dw^i_t}$, for a Brownian motion $w$. As the probability measure $\PP^{x,y}_\ep$ has finite energy \cite{AiraultMalliavin}, a theorem of Sugita, theorem 4.2 in \cite{Sugita}, ensures that we have $\big\{\PP^{x,y}_\ep(A)\big\}^p \leq c\,\textrm{\emph{C}}_p^r(X^\ep\in A)$, for some positive constant $c$ and all Borel sets $A$ in Wiener space; so a large deviation upper bound for $\textrm{\emph{C}}_p^r$ implies a corresponding result for $\PP^{x,y}_\ep(\cdot)$. It does not seem possible to get the large deviation lower bound by these methods.  \vspace{0.2cm} 
   \item We shall see in section \ref{SectionFirstOrderAsymptotics} that the large deviation principle stated in theorem \ref{MainLDPThm} leads directly to the identification of the first order asymptotics of $\PP^{x,y}_\ep$ under some mild conditions on $(x,y)$, in the sense that $\PP^{x,y}_\ep$ converges weakly to a Dirac mass on some particular path $\gamma$ from $x$ to $y$. It is natural in that setting to push further the analysis and try and get a second order asymptotics. This is done in the forthcoming work \cite{BailleulMesnagerNorris} where it is proved that the fluctuation process around the deterministic limit $\gamma$ is a Gaussian process whose covariance involves the (non-constant rank) sub-Riemannian geometry associated with the operator $\mcL$. This requires that the pair $(x,y)$ lies outside some intrinsic cutlocus associated with $\mcL$.
\end{enumerate}
\end{rems}

\section{Proof of the large deviation principle}
\label{SectionProofLDP}

The proof of theorem \ref{MainLDPThm} follows the pattern of proof devised by Hsu in \cite{Hsu} to prove a similar large deviation principle in a Riemannian setting where $\mcL$ is the Laplacian of some Riemannian metric on $M$. Our reasoning relies crucially on L{\'e}andre's logarithmic estimate \cite{Leandre1}, \cite{Leandre2}
\begin{equation}
\label{EqLeandreEstimate}
\underset{\ep\searrow 0}{\lim}\;\epsilon \log p_\epsilon(z,z') = -\frac{d^2(z,z')}{2},
\end{equation}
which holds uniformly \wrt $(z,z')\in M^2$, as well as on Sanchez-Calle's estimate 
\begin{equation}
\label{EqSanchezCalle}
p_t(z,z')\leq c\,t^{-m}, 
\end{equation}
which holds for some positive constant $c$ and all $z,z'\in M$ and $t>0$, see \cite{Sanchez}. 

\ssk

Write $\Omega^x$ for the set of continuous paths $\omega : [0,1]\rightarrow M$ started from $x$; we equip $\Omega^x$ and $\Omega^{x,y}$ with the metric of uniform convergence inherited from the ambiant space. Fix $\alpha\in\big(\frac{1}{3},\frac{1}{2}\big)$.

\bigskip

\noindent \textbf{a) Exponential tightness of the family of probability measures $(\PP^{x,y}_\epsilon)_{0<\epsilon\leq 1}$ on $\mcC^\alpha_{x,y}\big([0,1],M\big)$.}  Given $n=n(N)\geq 7$ and $K=K(N)$, to be fixed later as functions of some parameter $N$, we define a compact subset $C_N$ both of $\Omega^{x,y}$ and $\mcC^\alpha_{x,y}\big([0,1],M\big)$ setting
$$
C_N = \Big\{\omega\in\Omega^{x,y}\,;\,\underset{0<t-s\leq \frac{1}{n}}{\sup}\,\frac{|\omega_t-\omega_s|_d}{|t-s|^\alpha}\leq K\Big\}.
$$
The above supremum is over the set of all times $s,t\in [0,1]$. We first work on the time interval $[0,2/3]$ to evaluate the $\PP^{x,y}_\ep$-probability of $C_N$, to avoid the difficulties coming from the singularities of the drift at time $1$, in the classical dynamical description of the bridge as the solution to a stochastic differential equation. Set
$$
(\star) := \PP^{x,y}_\epsilon\left(\underset{s,t\in [0,2/3],\;0<t-s\leq \frac{1}{n}}{\sup}\,\frac{|\omega_t-\omega_s|_d}{|t-s|^\alpha}> K\right) \leq \frac{n}{2}\,\underset{0\leq r\leq 2/3}{\sup}\,\PP^{x,y}_\epsilon\left(\underset{r\leq s<t\leq r+2/n}{\sup}\,\frac{|\omega_t-\omega_s|_d}{|t-s|^\alpha}> K\right).
$$
Using \eqref{EqConditonnedProbability} and the Markov property provides the upper bound 
\begin{equation}
\begin{split}
\label{EqUpperBoundProbOscillations}
(\star) &\leq \underset{0\leq r\leq \frac{2\ep}{3}}{\sup}\,\EE^x\left[\frac{p_{\epsilon-r-\frac{2\ep}{n})}(\omega_{\epsilon(r+\frac{2\ep}{n})},y)}{p_\epsilon(x,y)}\,;\,\underset{r\leq s <t\leq r+\frac{2\epsilon}{n}}{\sup}\frac{|\omega_t-\omega_s|_d}{|t-s|^\alpha}\geq K\right] \\
&\leq \;\frac{c\epsilon^{-m}}{p_\epsilon(x,y)}\,\underset{z\in M}{\sup}\,\PP^z\Big(\underset{0\leq s <t\leq \frac{2\epsilon}{n}}{\sup}\frac{|\omega_t-\omega_s|_d}{|t-s|^\alpha}\geq K\Big).
\end{split}
\end{equation}
By Lyons' universal limit theorem, as stated for instance under the form given in theorem 11 in \cite{RoughFlows}, there exists universal controls on the oscillation of solutions of stochastic differential equations in terms of the oscillations of Brownian motion and its L{\'e}vy area. More precisely, there exists positive constants $a_i, b_i$, depending only on the vector fields $V,V_i$, such that 
\begin{equation}
\begin{split}
\underset{z\in M}{\sup}\,\PP^z\Big(\underset{0\leq s <t\leq \frac{2\epsilon}{n}}{\sup}\frac{|\omega_t-\omega_s|}{|t-s|^\alpha}\geq K\Big) &\leq a_1\Big\{\textrm{{\bf P}}\Big(\big\|{\bf B}_{[0,(2\epsilon)/n]}\big\|\geq b_1K\Big) + \textrm{{\bf P}}\Big(\big\|{\bf B}_{[0,(2\epsilon)/n]}\big\|^3\geq K\wedge\frac{n}{3}\Big)\Big\}, \\
&\leq a_2 \textrm{{\bf P}}\Big(\big\|{\bf B}_{[0,(2\epsilon)/n]}\big\| \geq b_2(K\wedge n)^{1/3}\Big)
\end{split}
\end{equation}
where ${\bf B}_{[0,(2\epsilon)/n]}$ is the Brownian $\frac{1}{\alpha}$-rough path on the time interval $\big[0,\frac{2\epsilon}{n}\big]$, defined on some probablity space $({\bf \Omega},\mcF,\textrm{{\bf P}})$, and $\big\|{\bf B}_{[0,(2\epsilon)/n]}\big\|$ stands for the homogeneous rough path norm of ${\bf B}_{[0,(2\epsilon)/n]}$; see for instance chapter 10.1 of \cite{FVBook}. It follows from the equality in law $\big\|{\bf B}_{[0,(2\epsilon)/n]}\big\| = \sqrt{\frac{2\ep}{n}}\big\|{\bf B}_{[0,1]}\big\|$, the Gaussian character of $\big\|{\bf B}_{[0,1]}\big\|$ under $\textrm{\bf P}$, and L{\'e}andre's estimate \eqref{EqLeandreEstimate} for $p_\epsilon(x,y)$, that  
$$
\epsilon\log\PP^{x,y}_\epsilon\left(\underset{s,t\in [0,2/3],\;0<t-s\leq \frac{1}{n}}{\sup}\,\frac{|\omega_t-\omega_s|}{|t-s|^\alpha}> K\right) \leq \frac{d^2(x,y)}{2} + o_\ep(1) -\frac{n(K\wedge n)^{2/3}}{2}\,b_2^2,
$$
so we have
\begin{equation}
\label{EqALmostTight}
\overline{\underset{\epsilon\searrow 0}{\lim}} \;\epsilon\log\PP^{x,y}_\epsilon\left(\underset{s,t\in [0,2/3],\;0<t-s\leq \frac{1}{n}}{\sup}\,\frac{|\omega_t-\omega_s|}{|t-s|^\alpha}> K\right) \leq -N
\end{equation}
by choosing $n=n(N)$ and $K=K(N)$ big enough. 

\ssk

To get a similar estimate when working on the whole time interval $[0,1]$, remark that since $M$ is compact and the operator $\mcL$ is hypoelliptic, it has a smooth positive invariant measure. If we use this measure as our reference measure $\Vol$, then $\widehat{p}_t(z,z') = p_t(z',z)$ is the heat kernel of another operator $\widehat{\mcL}$ which satisfies the same conditions as $\mcL$. Write $\widehat{\PP}^{z,z'}_\ep$ for the law of the associated bridge. So the class of measures $\Big(\PP^{z,z'}_\epsilon\Big)_{z\neq z' \in M}$ constructed from hypoelliptic operators $\mcL$ as in \eqref{EqDefnL}, satisfying the Lie bracket assumption, is preserved under time reversal. Applying inequality \eqref{EqALmostTight} to the measure $\widehat{\PP}^{y,x}_\ep$ on $\Omega^{y,x}$ obtained by time-reversal of $\PP^{x,y}_\ep$, we conclude with \eqref{EqALmostTight} that 
$$
\overline{\underset{\epsilon\searrow 0}{\lim}} \;\epsilon\log\PP_\epsilon^{x,y}(C_N^c)\leq -N.
$$
So the family $\big(\PP_\epsilon^{x,y}\big)_{0<\ep<1}$ of probabilities on $\mcC_\al^{x,y}\big([0,1],M\big)$ is exponentially tight, which proves in particular point 1. As the inclusion of $\mcC_\al^{x,y}\big([0,1],M\big)$ into $\big(\Omega^{x,y},\|\cdot\|_\infty\big)$ is continuous, it suffices, by the inverse contraction principle, to prove that $\big(\PP_\epsilon^{x,y}\big)_{0<\ep<1}$ satisfies a large deviation principle in $\big(\Omega^{x,y},\|\cdot\|_\infty\big)$, with good rate function $J$, to prove point 2 of the theorem, in so far as $J$ is also a good rate function on $\mcC_\al^{x,y}\big([0,1],M\big)$. We follow closely Hsu's work \cite{Hsu} to prove that fact.

\bigskip

\noindent \textbf{b) Large deviation upper bound for $(\PP^{x,y}_\epsilon)_{0<\epsilon\leq 1}$.} We first prove the upper bound for a compact subset $C$ of $\Omega^{x,y}$. For $0<a<1$, set 
$$
C^a = \big\{\omega\in\Omega^{x,y}\,;\,\exists\,\rho\in C\textrm{ such that } \omega_s = \rho_s, \textrm{ for }0\leq s\leq 1-a\big\}
$$
and 
$$
C^a_* = \big\{\omega\in\Omega^x\,;\,\exists\,\rho\in C\textrm{ such that } \omega_s = \rho_{(1-a)s}, \textrm{ for all }0\leq s\leq 1\big\}.
$$
The set $C^a$ is closed in both $\Omega^x$ and $\Omega^{x,y}$, and $C\subset C^a$. Using \eqref{EqConditonnedProbability} and the Markov property, we get as in \eqref{EqUpperBoundProbOscillations} the inequality
\begin{equation*}
\begin{split}
\PP^{x,y}_\epsilon(C) &\leq \PP^{x,y}_\epsilon(C^a) \leq \EE^x_\ep\left[\frac{p_{a\epsilon}(\omega_1,y)}{p_\epsilon(x,y)}{\bf 1}_{\omega\in C^a_*}\right] \\
&\leq \frac{c\epsilon^{-m}}{p_\epsilon(x,y)} \PP^x_\ep(C^a_*).
\end{split}
\end{equation*}
As $C^a_*$ is closed in $\Omega^x$, we have by the classical Freidlin-Wentzell large deviation principle for $\PP^x_\ep$
$$
\underset{\epsilon\searrow 0}{\limsup}\,\epsilon\log \PP^{x,y}_\epsilon(C) \leq \frac{d^2(x,y)}{2}-\frac{1}{1-a}\,\underset{\omega\in C^a_*}{\inf}\,I(\omega).
$$
Using the lower semicontinuity of $I$ on $\Omega^x$, it is straightforward to use the compacity of $C$ to see that $\underset{a\searrow 0}{\limsup}\,\underset{\omega\in C^a_*}{\inf}\,I(\omega)\geq \underset{\omega\in C}{\inf}\,I(\omega)$, as done in \cite{Hsu}, p.112. This proves the upper bound 
$$
\underset{\epsilon\searrow 0}{\limsup}\,\epsilon\log \PP^{x,y}_\epsilon(C) \leq -\inf_C J,
$$
for a compact set $C$; it is classical that the exponential tightness proved in point a) implies in that case the upper bound for any closed set.

\bigskip

\noindent \textbf{c) Large deviation lower bound for $(\PP^{x,y}_\epsilon)_{0<\epsilon\leq 1}$.} We use the notation $\|f\|_{[a,b]}$ to denote the uniform norm of some function $f$ defined on some time interval $[a,b]$. Given an open set $U$ in $\Omega^{x,y}$,we aim at proving that we have 
\begin{equation}
\label{EqLDPLowerBounder}
\underset{\ep\searrow 0}{\liminf}\,\ep\log\PP^{x,y}_\ep(U) \geq -J(\gamma)
\end{equation}
for any $\gamma\in U$ with finite energy $I(\gamma)$. Pick such a path $\gamma\in U$ and $b>0$ small enough for the ball in $\Omega^x$ with center $\gamma$ and radius $b$ to be included in $U$. Set for $0<a<1$
\begin{equation*}
U^{a,b} = \big\{\omega\in\Omega^{x,y}\,;\,\|\omega-\gamma\|_{[0,1-a]} < b\big\},\quad F^{a,b} = \big\{\omega\in\Omega^{x,y}\,;\,\|\omega-\gamma\|_{[1-a,1]} \geq b\big\}
\end{equation*}
and $U^{a,b}_* = \big\{\omega_*\in\Omega^x\,;\,\exists\,\omega\in U\textrm{ such that } \omega_*(s) = \omega_{(1-a)s}, \textrm{ for all }0\leq s\leq 1\big\}$.  We have $U^{a,b}\subset \big(U\cup F^{a,b}\big)$, so $\PP^{x,y}_\ep(U) \geq \PP^{x,y}_\ep\big(U^{a,b}\big) - \PP^{x,y}_\ep\big(F^{a,b}\big)$. We prove \eqref{EqLDPLowerBounder} by showing that $\underset{\ep\searrow 0}{\liminf}\,\ep\log\PP^{x,y}_\ep\big(U^{a,b}\big) \geq -J(\gamma)$, and $\underset{\ep\searrow 0}{\liminf}\,\ep\log\PP^{x,y}_\ep\big(F^{a,b}\big) = -\infty$.

\ssk

Given $\lambda>0$, write $B_\lambda(y)$ for the sub-Riemannian open ball in $M$, with center $y$ and radius $\lambda$. Using the Markov property as above, we have
\begin{equation*}
\begin{split}
\PP^{x,y}_\ep\big(U^{a,b}\big) &= \EE^{x,y}_\ep\Big[\PP^{x,X_{1-a}}_{\ep(1-a)}\big(U^{a,b}_*\big)\Big] \geq \int \PP^{x,z}_{\ep(1-a)}\big(U^{a,b}_*\big)\frac{p_{\ep(1-a)}(x,z)p_{\ep a}(z,y)}{p_\ep(x,y)}{\bf 1}_{z\in B_\lambda(y)}\,dz \\
&\geq \frac{\underset{z\in B_\lambda(y)}{\min}p_{\ep a}(z,y)}{p_\ep(x,y)} \int \PP^{x,z}_{\ep(1-a)}\big(U^{a,b}_*\big) {\bf 1}_{z\in B_\lambda(y)}p_{\ep(1-a)}(x,z)\,dz \\
&\geq \frac{\underset{z\in B_\lambda(y)}{\min}p_{\ep a}(z,y)}{p_\ep(x,y)}\,\PP^{x}_{\ep (1-a)} \big(U^{a,b}_*\cap\{\omega_1\in B_\lambda(y)\}\big)
\end{split}
\end{equation*}
Define $\gamma_a(s) = \gamma_{(1-a)s}$ for all $0\leq s\leq 1$. As $\gamma$ has finite energy, one can pick some control $h\in H^1_0$ such that $\gamma^h=\gamma$; we have $d(\gamma_a(1),y) \leq \int_{1-a}^1|\dot h_s|_\ell\,ds \leq \sqrt{a\int_{1-a}^1|\dot h_s|^2_\ell\,ds}$. The choice of $\lambda = \lambda(a) = 2\sqrt{a\int_{1-a}^1|\dot h_s|^2_\ell\,ds}$ ensures that the open set $U^{a,b}_*\cap\big\{\omega_1\in B_\lambda(y)\big\}$ contains $\gamma_a$, so it is nonempty; also, $\frac{\lambda(a)^2}{a}\rightarrow 0$ as $a$ tends to $0$. Using the classical Freidlin-Wentzell large deviation theory and the uniform character of L\'eandre's estimate \eqref{EqLeandreEstimate}, the above lower bound for $\PP^{x,y}_\ep\big(U^{a,b}\big)$ gives
$$
\underset{\ep\searrow 0}{\liminf}\,\ep\log\PP^{x,y}_\ep\big(U^{a,b}\big) \geq \frac{-I(\gamma_a)}{1-a}+\frac{d(x,y)^2}{2} - \frac{\lambda(a)^2}{2a}, 
$$ 
from which the inequality $\underset{\ep\searrow 0}{\liminf}\,\ep\log\PP^{x,y}_\ep\big(U^{a,b}\big) \geq -J(\gamma)$ follows, since $I(\gamma_a)\rightarrow I(\gamma)$ and $\frac{\lambda(a)^2}{a}\rightarrow 0$ as $a$ tends to $0$. 

\ssk

We now deal with the term $\PP^{x,y}_\ep\big(F^{a,b}\big)$. Set $\overline\gamma_s = \gamma_{1-s}$, for $0\leq s\leq 1$, and choose $a$ small enough to have $\big\|\overline\gamma-y\big\|_{[0,a]} \leq \frac{b}{2}$. We use the same time reversal trick and notations as above to estimate $\PP^{x,y}_\ep\big(F^{a,b}\big)$. Write
\begin{equation*}
\begin{split}
\PP^{x,y}_\ep\big(F^{a,b}\big) &= \widehat\PP^{y,x}_\ep\big(\big\|\omega-\overline\gamma\big\|_{[0,a]}\geq b\big) \leq \widehat\PP^{y,x}_\ep\big(\big\|\omega-y\big\|_{[0,a]}\geq \frac{b}{2}\big) \\
&\leq \frac{c\ep^{-m}}{p_\ep(y,x)} \widehat\PP^y_\ep\big(\big\|\omega-y\big\|_{[0,a]}\geq \frac{b}{2}\big).
\end{split}
\end{equation*}
L\'eandre's estimate \eqref{EqLeandreEstimate} and the classical large deviation results for $\widehat\PP^y_\ep$ give the existence of a positive constant $c$ such that we have
$$
\underset{\ep\searrow 0}{\liminf}\,\ep\log\PP^{x,y}_\ep\big(F^{a,b}\big) \leq \frac{d(x,y)^2}{2} - \frac{c}{a};
$$
this upper bound tends to $-\infty$ as $a$ tends to $0$. Points a), b) and c) all together prove theorem \ref{MainLDPThm}.

\bigskip

\section{First order asymptotics for bridges of degenerate diffusion processes}
\label{SectionFirstOrderAsymptotics}

Theorem \ref{MainLDPThm} provides a straightforward mean for investigating the first order asymptotics of $\PP^{x,y}_\ep$ as $\ep$ tends to $0$, for $x$ and $y$ in generic positions.

\begin{thm}[First order asymptotics of $\PP^{x,y}_\ep$]
If there exists a unique path $\gamma$ with minimal energy from $x$ to $y$, then $\PP^{x,y}_\ep$ converges weakly in $\big(\Omega^{x,y},\|\cdot\|_\infty\big)$ to a Dirac mass on $\gamma$ as $\ep$ tends to $0$.
\end{thm}

\begin{Dem}
We follow the proof of lemma 3.1 in \cite{Hsu}. Since the family $(\PP^{x,y}_\epsilon)_{0<\epsilon\leq 1}$ is tight by point a) in section \ref{SectionProofLDP}, let $\QQ$ be any limit measure. Given $b>0$, set
$$
C^b_N = C_N \cap\{\omega\in\Omega^{x,y}\,;\,\|\omega-\gamma\|_\infty>b\};
$$
then $\underset{\omega\in C^b_N}{\inf} \, J(\omega)>0$. Indeed, since the paths of $C^b_N$ are equicontinuous, if the infimum were null, we could extract from any sequence of paths $(\omega_n)_{n\geq 0}$ such that $J(\omega_n)$ converges to $0$ a uniformly converging subsequence with limit $\omega\in \overline C^b_N$ say. We should then have $J(\omega)=0$, by the lower semicontinuity of $J$, that is $\omega=\gamma$, since there is a unique path from $x$ to $y$ with minimal energy, in contradiction with the fact that elements of $\overline C^b_N$ satisfy the inequality $\|\omega-\gamma\|_\infty\geq b>0$. As a consequence, the above large deviation upper bound implies 
$$
\QQ\big(C^b_N\big)\leq \underset{\ep\searrow}{\liminf}\;\PP^{x,y}_\ep\big(\overline C^b_N\big) = 0;
$$
sending $N$ tend to infinity, it follows that 
$$
\QQ\big(\omega\in\Omega^{x,y}\,;\,\|\omega-\gamma\|_\infty>b\big) = 0.
$$
As this holds for all $b>0$, we have $\QQ = \delta_\gamma$, from which the convergence of $\PP^{x,y}_\epsilon$ to $\delta_\gamma$ follows.
\end{Dem}

\ssk

\noindent Note that the set of pairs of points $(x,y)\in M^2$ such that $x$ and $y$ are joined by a unique path of minimal energy is dense in $M^2$.

\end{document}